\documentclass[10pt]{article}
\usepackage{amsmath,amssymb,amsthm,amscd}
\numberwithin{equation}{section}

\def\p{\partial}
\def\b{\bar}

\def\i{\sqrt{-1}}
\def\vphi{\varphi}

\def\o{\omega}

\def\cC{{\cal C}}

\def\cE{{\cal E}}

\def\cH{{\cal H}}

\def\h{{\mathfrak h}}

\def\cC{{\mathcal C}}

\def\cE{{\mathcal E}}

\def\cG{{\mathcal G}}
\def\cH{{\mathcal H}}

\newtheorem{prop}{Proposition}[section]
\newtheorem{theo}[prop]{Theorem}

\newtheorem{cor}[prop]{Corollary}
\newtheorem{rem}[prop]{Remark}

\newtheorem{defi}[prop]{Definition}

\newtheorem{claim}[prop]{Claim}

\def\and{\quad{\rm and}\quad}

\let\lra=\longrightarrow

\def\mapright\#1{\,\smash{\mathop{\lra}\limits^{\#1}}\,}

\begin {document}
\bibliographystyle{plain}
\title{On the Calabi flow}
\author{X. X. Chen\footnote{The author is partially supported by a NSF grant. } \;
and  W. Y. He\footnote{The author is partially supported by a NSF supplement grant.}\\ Department of Mathematics\\
University of Wisconsin,  Madison}
\date{}
\maketitle
\tableofcontents

\section{Introduction}

~~~~In \cite{Ca01}, E. Calabi studied the variational problem of
minimizing the so-called ``Calabi energy" \footnote{It is the
$L^2$ norm of the scalar curvature function of the K\"ahler
metric.} in any fixed cohomology class of K\"ahler metrics. Any
smooth critical point is called either constant scalar curvature
(CscK) metric or extremal K\"ahler (CextK) metric depending on
whether the Futaki Character vanishes or not in this class. CextK
metrics include the more famous K\"ahler Einstein metric as a
special case when the K\"ahler class is the canonical K\"ahler
class. In recent years, the study of extremal K\"ahler (or CscK)
metric has attracted intensive attention and many important works
emerge. In particular, the uniqueness problem is completely
settled, \cite{chen001}, \cite{Dona01}, \cite{Ma01}, and
\cite{chentian005}. For existence, the situation is much more
complicated. In 1976, S. T. Yau solved the famous Calabi
conjecture which implies that any K\"ahler manifold with vanishing
first Chern class has a Calabi-Yau metric, that is, a Ricci-flat
K\"ahler metric. Around the same time, T. Aubin and S. T. Yau
independently proved existence of K\"ahler-Einstein metrics on
compact K\"ahler manifolds with negative first Chern class.  In
\cite{Tian90},  G. Tian proved that a complex surface with
positive first Chern class admits a K\"ahler-Einstein metric if
and only if its automorphism group is reductive. In 1997, a  new
algebraic invariant (K stability) was introduced by G. Tian
\cite{Tian97} as an obstruction to the existence of K\"ahler
Einstein metric.  Later, S. K. Donaldson \cite{Dona01} gave a new
definition of K-stability by using weights of Hilbert points...
These invariants can also be extended to be the obstruction to the
existence of Calabi's Extremal K\"ahler (CextK) metrics or CscK
metrics, c.f. \cite{Ma02} and \cite{PT} for further references.\\

In 1984, Calabi constructed the first CextK metric, in
$\mathbb{CP}^n$ blown up at a point, by solving a 4th order ODE
(the so-called ``Calabi's anstaz").  There are many beautiful
works in constructing special CscK metric or CextK metric on many
K\"ahler manifolds (cf. Lebrun-Simanca  \cite{Lebrunsmanca}
\cite{Lebrunsinger} \cite{AposCalGaudFrid}  and
\cite{Rollinsinger} etc and reference therein).  More recently, J.
Fine \cite{Fine1} has used adabiatic limit techniques to construct
non-trivial CscK metrics on complex surfaces.  Around the same
time, using method of gluing and grafting Arezzo-Pacard
\cite{Arezzopacard1}\cite{Arezzopacard2} were able to construct
new CscK metrics on the blow up of manifold if the original
manifold has a CscK metric. In \cite{DonaFine}, Donaldson-Fine
classify the toric anti-self-dual 4-manifold with some simple
holomorphic data by using twistor correspondence.\\

However, there has not been much progress made on the existence of
general Calabi's extremal K\"ahler metrics, via direct PDE method.
Adopting the strategy of direct variational approach, the first
author \cite{chenthesis} studied the problem of minimizing the
Calabi energy on Riemann surface, in an attempt to understand how
compactness fails in general. Recent work of S. Donaldson
\cite{Dona02} on toric surfaces, while still a special case, leads
hope that the existence problem might be approachable on toric
surfaces. Another ambitious program, not directly K\"ahler but
very much related, is the recent work of Cheeger-Tian \cite{ct01}
and Tian-Viaclosky \cite{tv01}.  However, for general K\"ahler
manifold without boundary, the existence problem is still too
difficult. Partially it is because the equation for CscK metric is
a fully nonlinear 4th order partial differential equation.
Envisioning this difficulty, in the same paper defining CscK
metrics, E. Calabi proposed the so-called Calabi flow to attack
the existence problem.  Unlike the K\"ahler Ricci flow, (in Calabi
flow) one deforms the K\"ahler potential in the direction of the
scalar curvature.  In \cite{Ca01}, \cite{Ca02}, Calabi shows that
a) The Calabi energy is decreasing along the Calabi flow; b) The
Calabi energy is weakly convex near a CscK or CextK metric.
Therefore, at least conceptually, there is a good chance that the
Calabi flow will converge in some fashion. These two properties
play a crucial role in understanding the stability problem of the
Calabi flow.

In general dimensional case, very few results are known or written
down, even though several authors have comments on the importance
of such a flow (cf. \cite{Dona03}). Even for the short time
existence, which is known to the experts (because the flow is
parabolic) for smooth initial data, the precise/optimal result was never formulated in
literature. In this paper, we first give a short time existence
result:

\begin{theo}Let $(M, [\omega])$ be a polarized K\"ahler manifold without boundary.
For any K\"ahler potential with $C^{3,\alpha}(M, \omega_g)$ norm
bounded and such that the corresponding  metric is uniform
equivalent to $\omega_g$, then the Calabi flow exists for a short
time and the K\"ahler metric becomes immediately smooth for $t >
0.\;$\end{theo}


In finite dimensions, once one establishes a contracting process
to derive the short time existence, then it is immediately clear
that, if one starts near a fixed point, it will stay in small
neighborhood of the fixed point and eventually converge to the
fixed point fast.  However, in infinite dimensional case, the
picture is much worse and peculiar phenomenon could occur.  For
certain parabolic flow with short time existence property, no
matter how close you start from a fixed point, once the flow
starts, it will diverge quickly away from fixed point. The hope
for establishing a Stability theorem lies in the fact that the
Hessian of the Calabi energy is strictly positive (unless it is in
the direction of holomorphic isometry).

\begin{theo} (Stability theorem)
Suppose $g$ is a CscK metric in $[\omega]$ on $M$. If the ${C^{3,
\alpha}(M)}$ norm of the initial potential (w.r.t $g$) is small
enough, then the flow will exist for all time and
$g(t)=g_{\varphi(t)}$ will converge to a CscK metric $g_\infty$ in
the same class $[\o]$ in $C^{\infty}$ sense. Moreover, the Calabi
energy will decay exponentially fast and consequently the
convergence is exponentially fast along the flow.
\end{theo}

\begin{rem}The limit CscK metric $g_\infty$ could be different from
the original CscK metric $g$ in holomorphic coordinates. By the
recent result of \cite{chentian005}, the CextK metric is unique in
the fixed K\"ahler class up to automorphism. And so $g$ and
$g_\infty$ differ by an automorphism. See the example in section 4
for more details.
\end{rem}

To attack the long time existence of the Calabi flow, in the early
90s, the first author concentrated on the weak compactness of
K\"ahler metrics with uniform bound on the Calabi energy
\cite{chen002}, \cite{chen003} and \cite{chen004}. This is a
crucial step in understanding the whole picture of the Calabi flow
since the Calabi energy is decreasing along the flow.  Using this
weak compactness theorem, the first author \cite{Chen} was able to
prove the flow exists for all time in Riemann surface and the flow
will converge to a CscK metric exponentially fast assuming the
uniformization theorem.  It gives a new proof to a theorem of P.
Chru\'sciel \cite{Ch} on the Calabi flow in Riemann surface from
a completely different perspective.\\

In this paper, we follow the path of \cite{Chen} to tackle the
long time existence under some suitable curvature assumption. In a
fixed K\"ahler class, we obtain a compactness theorem under the
uniform bound of Ricci curvature and potential.

\begin{theo}
(Compactness theorem) All metrics
$\omega_\varphi=\omega+\sqrt{-1}\p \bar{\p }\varphi$ in the space
of K\"ahler metrics with both the potential $\varphi$ and the
Ricci curvature $Ric_\varphi$ uniformly bounded are equivalent and
compact in $C^{1, \alpha}-$topology for any $\alpha \in (0,1)$ (It
is equivalent to say $\varphi$ is uniformly bounded in $C^{3,
\alpha}$ for any $\alpha \in (0,1)$).
\end{theo}

As a consequence, we prove that the Calabi flow  will exist for
all time if the evolving Ricci curvature is uniformly bounded.

\begin{theo}For the Calabi flow initiating from any smooth K\"ahler metric,
the flow exists as long as the Ricci curvature stays uniformly
bounded.
\end{theo}

\begin{rem}For Hamilton's Ricci flow \cite{Hamilton82},  it was
proved by Hamilton that the flow will continue as long as the Riemannian
curvature stays uniformly bounded. After Perelman's work  \cite{perelman01}, \cite{perelman02},
 N. Sesum \cite{Se01} was able to improve this important result saying that the Ricci
flow will continue as long as the Ricci curvature is bounded.
Readers are referred to \cite{CLP}  for a complete updated reference on this important
subject.
\end{rem}
In \cite{Chen}, the weak compactness of the space of K\"ahler
metrics with uniform  Calabi energy and area in Riemann surface is
used critically in obtaining the long time existence of the Calabi
flow. Motivated by \cite{Chen}, 
it is highly desirable to obtain a weak compactness theorem
similar to the corresponding theorems in \cite{chen002},
\cite{chen003}, and \cite{chen004}  in general
K\"ahler manifold.  This will be discussed in a subsequent paper.\\

In \cite{CC}, Calabi-Chen showed that the Calabi flow essentially
decreases the geodesic distance in the space of K\"ahler metrics.
This property must play an important role in the future study of
the Calabi flow.  In particular, this property suggests that the
flow shall exist globally (cf. \cite{chen05}) and the flow will
converge to a CextK metric except a codimension 2 subvariety
although the complex structure of the limiting manifold might be
different from the original one.  In \cite{Dona03}, Donaldson
describes precisely what the limit of the Calabi flow shall be
under various situation, from perspective of sympletic moment map.
Motivated by these discussions, it is very important to study the
so-called ``removing singularity" for CscK metric or CextK metric.
\begin{theo} Let $g$ be a weakly CscK metric in a punctured disc which is smooth in a small neighborhood
of the boundary of the disc.  If $g$ is uniformly bounded from
above and below (with respect to the Euclidean metric) in the
punctured disc and its tensors are weakly in $W^{1,2}_{loc}$ ,
then it is a smooth CscK K\"ahler metric in the entire disc.
\end{theo}
\begin{rem}We simplify the proof by requiring the boundary data to be smooth in a small
neighborhood of the boundary of the disc and provide a complete
proof the last claim (Claim 6.6). Since the aim of this theorem is
about removing possible singularities, the significance of the
theorem is not affected.  On the other hand, we add a remark that
to prove the full version (we mean that the boundary data is
smooth instead of smooth in a small neighborhood of the boundary),
one need add one more step which is about to solve the
Monge-Ampere equation \[ \log{\det{(\delta_{kl}+\varphi_{k\bar
l})}}=f\] in a strongly pseudo-convex domain with $f\in C^\alpha$.
\end{rem}

This is a first step in establishing a more ambitious regularity
result (cf \cite{chen05}) for weak CscK metrics on K\"ahler
manifolds where the ``suspected" singular locus is a union of
subvarieties of codimension $2$ or higher.  In this regard, the
assumption of uniform ellipticity is necessary in the following
sense:  if we blow up $M$ at a point $p$, we can denote the
resulting K\"ahler manifold as $\hat M$ and the exceptional
divisor  $\cE.\;$ Suppose $M$ and $\hat M$ is each equipped with a
CscK metric.  Then, one can view either of these CscK metrics as a
weak CscK metric in the wrong manifold and singularities cannot be
removed. One easily observes in both cases, the uniformly
ellipticity is violated, but in different directions. In other
words, one can perhaps weaken  the assumption of uniform
ellipticity in some way. However, some form of ellipticity must be presented.\\

\noindent {\bf Acknowledgement} The first author is grateful for
the brief, but insightful conversation with S. K. Donaldson on the
issue of Stability of the Calabi flow. The first author also
benefits from discussion of removing singularity for CscK metrics
with F. Pacard.  He wishes to thank both of them. The second
author would like to thank S. B. Angenent and P. Rabinowitz for
the conversation of the local existence and stability of general
parabolic equation. In particular, the elliptic operator
$A(\varphi)$ in Section 3 is suggested by S. B. Angenent.  Both
authors are grateful for X. H. Zhu's kind comments on an earlier
version of this paper.

\section{Preliminary}
\subsection{Notations in K\"ahler Geometry}
~~~~Let $M$ be a compact complex manifold of complex dimension
$n$. An Hermitian metric $g$ on $M$ in local coordinates is given
by \begin{eqnarray*} g=g_{i\bar{j}}dz^i\otimes
dz^{\bar{j}},\end{eqnarray*} where $\{g_{i\bar{j}}\}$ is a
positive definite Hermitian matrix smooth function. And we use
$\{g^{i\bar{j}}\}$ to denote the inverse matrix of
$\{g_{i\bar{j}}\}$. The K\"ahler condition says that the
corresponding K\"ahler form $\omega=\i g_{i\bar{j}}dz^i\wedge
dz^{\bar{j}}$ is a closed $(1,1)$ form. The K\"ahler class of
$\omega$ is its cohomology class $[\omega]$ in $H^2(M,
\mathbb{R})$. By the Hodge theory, any other K\"ahler form in the
same class is of the form
\begin{eqnarray*}\omega_{\varphi}=\omega+\i
\p\bar{\p}\varphi>0,\end{eqnarray*} for some real valued function
$\varphi$ on $M$, where
\begin{eqnarray*}\p\bar{\p}\varphi=\sum^n_{i,j=1}\p_i\p_{\bar{j}}\varphi
dz^i\wedge dz^{\bar{j}}=\varphi,_{i\bar{j}}dz^i\wedge
dz^{\bar{j}}.\end{eqnarray*} The corresponding K\"ahler metric is
denoted by
$g_\varphi=\left(g_{i\bar{j}}+\varphi,_{i\bar{j}}\right)dz^i\otimes
dz^{\bar{j}}$, and we use $\{g^{i\bar{j}}_\varphi\}$ to denote the
inverse matrix of $\{g_{i\bar{j}}+\varphi,_{i\bar{j}}\}.$ For
simplicity, we use both $g$ and $\omega$ to denote the K\"ahler
metric. Define the space of K\"ahler potentials \begin{eqnarray*}
\cH_\omega=\left\{\varphi|\omega_\varphi=\omega+\sqrt{-1}\partial\bar{\partial}\varphi>0,~~
\varphi \in C^{\infty}(M) \right\},\end{eqnarray*} which is called the space of K\"ahler metrics and is the main object we are interested in.\\

Given a K\"ahler metric $\omega$, its volume form is \[
\omega^n=\frac{(\i)^n}{n!}\det{(g_{i\bar{j}})}dz^1\wedge
dz^{\bar{1}}\wedge\cdots \wedge dz^n\wedge dz^{\bar{n}}.\] The
Ricci curvature of $\omega$ is locally given by \[
R_{i\bar{j}}=-\p_i\p_{\bar{j}}\log {\det{(g_{k\bar{l}})}}.\] Its
Ricci form is of the form \begin{eqnarray*} Ric_{\omega}=\i
R_{i\bar{j}}dz^idz^{\bar{j}} =-\i \p_i \p_{\bar
j}\log{\det{(g_{k\bar{l}})}}.\end{eqnarray*} It is a real, closed
$(1,1)$ form. The cohomology class of Ricci form is the famous
first Chern class $c_1(M)$, independent of the metric.

\subsection{The Calabi flow}
~~~~Given a polarized compact K\"ahler manifold $(M, [\omega])$,
for any $\varphi \in \cH_\omega,$ E. Calabi introduced the Calabi
functional in \cite{Ca01}, \cite{Ca02}, \begin{eqnarray*} \cC
a(\omega_\varphi)=\int_MR_{\varphi}^2\omega^n_\varphi,\end{eqnarray*}
where $R_\varphi$ is the scalar curvature of $\omega_\varphi$.
Note that both the total volume \begin{eqnarray*}
V_\varphi=\int_M\omega_\varphi^n\end{eqnarray*} and the total
scalar curvature
\begin{eqnarray*} S_\varphi=\int_MR_\varphi\omega_\varphi^n\end{eqnarray*} remain
unchanged when $\varphi$ varies in $\cH_\omega$. As a consequence,
the average of the scalar
curvature\begin{eqnarray*}\underline{R}=\frac{S_\varphi}{V_\varphi}\end{eqnarray*}
is a constant depending only on the class $[\omega]$. Usually we
use the following modified Calabi energy \begin{eqnarray*}
\tilde{\cC}
a(\omega_\varphi)=\int_M\left(R_{\varphi}-\underline{R}\right)^2\omega^n_\varphi\end{eqnarray*}
to replace $\cC a(\omega_\varphi)$ since they only differ by a
topological constant. E. Calabi studied the variational problem to
minimize $\cC a(\omega_\varphi)$ in $\cH_\omega$. The critical
point turns out to be either CscK or CextK metric depending on
whether the Futaki character vanishes or not. The Futaki character
$f=f_\varphi: \h(M)\rightarrow \mathbb{C}$ is defined on the Lie
algebra $\h(M)$ of all holomorphic vector fields of $M$ as
follows, \begin{eqnarray*}
f_\varphi(X)=-\int_MX(F_\varphi)\omega_\varphi^n,\end{eqnarray*}
where $X \in \h(M)$ and $F_\varphi$ is a real valued function
defined by
\begin{eqnarray*} F_\varphi=G_\varphi (R_\varphi).\end{eqnarray*} $G_\varphi$
is the Hodge-Green integral operator, and $F_\varphi=G_\varphi
(R_\varphi)$ is equivalent to $\triangle_\varphi
F_\varphi=R_\varphi-\underline{R}$, where $\triangle_\varphi$ is
the Laplace operator of the metric $\omega_\varphi.$ In
\cite{Ca02}, E. Calabi showed that the Futaki character
$f=f_\varphi$ is invariant when $\varphi$ varies in $\cH_\omega$,
and the critical point of $\cC a(\omega_\varphi)$ is CscK metric
when $f$ vanishes, otherwise CextK metric when $f$ is not zero.
Moreover, the second variation form (Hessian form) of $\cC
a(\omega_\varphi)$ is semi-positive definite at CscK metrics or
CextK metrics with a finite dimensional kernel (corresponding to
holomorphic vector fields).

The existence of CscK metric (or CextK metric) seems intractable
at the first glance since the equation is a fully nonlinear 4th
order equation (6th order for CextK metric). In \cite{Ca01}, E.
Calabi proposed the so-called Calabi flow to approach the
existence problem. The Calabi flow is the gradient flow of the
Calabi functional, defined as the following parabolic equation
with respect to a real parameter $t\geq 0$, \begin{eqnarray*}
\frac{\p}{\p
t}g_{i\bar{j}}(t)=\p_i\p_{\bar{j}}R_{g(t)}.\end{eqnarray*} On the
potential level, the Calabi flow is of the form \begin{eqnarray*}
\frac{\p \varphi}{\p t}=R_\varphi-\underline{R}.\end{eqnarray*}
Under the Calabi flow, we
have\begin{eqnarray*}\frac{d}{dt}\int_M(R_\varphi-\underline{R})^2\omega^n_\varphi=-2\int_M\left(D_\varphi
R_\varphi, R_\varphi\right) \omega^n_\varphi,\end{eqnarray*} where
$D_\varphi$ is Lichn\'erowicz operator with respect to
$\omega_\varphi.$ Lichn\'erowicz operator  $D$ is defined by
\begin{eqnarray*} Df={f,_{\alpha\beta}}^{\alpha\beta},\end{eqnarray*} where
the covariant derivative is with respect to $\omega$. So the
Calabi energy is strictly decreasing along the flow unless
$\omega_\varphi$ is CextK or CscK metric.

\section{Short time existence}
~~~~A straightforward calculation shows that \begin{eqnarray*}
  R_\varphi&=&g^{i\bar{j}}_{\varphi}R_{i\bar{j}}(\varphi)\\
  &=& -g^{i\bar{j}}_{\varphi}\p_i\p_{\bar{j}}\bigl\{
    \log \det(g_{k\bar{l}}+\p_k\p_{\bar{l}}\varphi)
    \bigr\}\\
  &=& -g^{i\bar j}_{\varphi}\p_i\left\{
     g^{k\bar l}_{\varphi} \p_{\bar{j}} (g_{k\bar{l}}+\p_k\p_{\bar{l}}\varphi)
    \right\} \\
  &=& -g^{i\bar j}_{\varphi}\left\{
     g^{k\bar l}_{\varphi} \p_{i}\p_{\bar{j}} (g_{k\bar{l}}+\p_k\p_{\bar{l}}\varphi) +
     \p_{i}g^{k\bar l}_{\varphi} \p_{\bar{j}} (g_{k\bar{l}}+\p_k\p_{\bar{l}}\varphi)
    \right\} \\
  &=& -A(\nabla\varphi,\nabla^{2}\varphi)\varphi + f(\nabla
  \varphi,\nabla^{2}\varphi,\nabla^{3}\varphi),
\end{eqnarray*}
where
\begin{equation}\label{3-1}
  A(\varphi)w = A(\nabla\varphi,\nabla^{2}\varphi)w=g^{i\bar j}_{\varphi}g^{k\bar l}_{\varphi}
         \nabla_{i}\nabla_{\bar j}\nabla_{k}\nabla_{\bar l}w,
\end{equation}
$\nabla, \bar{\nabla}$ are both covariant derivatives with respect
to $\omega$ and $f(\cdots)$ is some function of the derivatives of
order $\leq 3$ of $\varphi$. $A(\varphi)$ is a strictly elliptic
operator with coefficients depending on the derivatives of order
$\leq 2$ of $\varphi$, so the Calabi flow is a standard 4th order
quasilinear parabolic equation. The quasilinear parabolic equation
is well studied by many authors, and the existence of solution is
established. In some sense the short time existence of the Calabi
flow is standard and well known. We would like to carry out the
details here for our purpose to tackle long time existence.\\

One approach for local existence and regularity theory of
quasilinear parabolic equation is to use the maximal regularity
theory of \cite{DP1}, which is based on the use of continuous
interpolation spaces. In \cite{An}, \cite{CS1}, the maximal
regularity was extended to functions which admit a prescribed
singularity at $t=0$. This extension allows us to take advantage
of the quasilinear equation to establish the smoothing property.
It turns out that this version of short time existence is very
useful for the problem of the long time existence of the flow.

\subsection{Function space and maximal regularity}
~~~~This section follows \cite{CS1} and we will state a main
theorem in \cite{CS1} which is used to prove the short time
existence in next section. In the following we assume $\theta \in
(0,1)$, and $E$ is a Banach space, $J=[0, T]$ for some $T>0$. We
consider functions defined on $\widetilde{J}=(0, T]$ with a
prescribed singularity at $0$. Set
\begin{eqnarray*} C_{1-\theta}(J, E):=\left\{u\in C(\widetilde{J},
E); [t\rightarrow t^{1-\theta}u] \in C(J, E), \lim_{t\rightarrow
0^+}t^{1-\theta}|u(t)|_E=0\right\},\end{eqnarray*}
\begin{eqnarray*}|u|_{C_{1-\theta}(J, E)}:=\sup_{t\in
\widetilde{J}}t^{1-\theta}|u|_E,\end{eqnarray*} and
\begin{eqnarray*} C^1_{1-\theta}(J, E):=\left\{u \in
C^1(\widetilde{J}, E); u, \dot{u} \in C_{1-\theta}(J,
E)\right\}.\end{eqnarray*} It is easy to verify $C_{1-\theta}(J,
E)$, equipped with the norm $\|.\|_{C_{1-\theta}(J, E)}$, is a
Banach space and
$C^1_{1-\theta}(J, E)$ is a subspace.\\

Let $E_1, E_0$ be two Banach space such that $E_1$ is continuously
embedded into $E_0$. The set of bounded linear operators from
$E_1$ to $E_0$ is denoted by $\cal L$$(E_1, E_0)$. Any operator $A
\in \cal L$$(E_1, E_0)$ can be considered as an unbounded in $E_0$
with domain $E_1$. It can happen that an $A \in \cal L$$(E_1,
E_0)$, seen as unbounded operator in $E_0$, is closed and
generates a strongly continuous semigroup denoted by
$\exp(-tA)(0\leq t< \infty)$. And if the associated semigroup
$\exp(-tA)$ is an analytic semigroup, we say $A \in \cal H$$(E_1,
E_0)$. Set \begin{eqnarray*} E_0(J)&&:=C_{1-\theta}(J,
E_0),\end{eqnarray*} and \begin{eqnarray*}
E_1(J)&&:=C^1_{1-\theta}(J, E_0)\cap C_{1-\theta}(J,
E_1),\end{eqnarray*} where
\begin{eqnarray*}|u|_{E_1(J)}:=\sup_{t\in
\widetilde{J}}t^{1-\theta}\left(|\dot{u}|_{E_0}+|u|_{E_1}\right).\end{eqnarray*}
By using $E_1(J), E_0(J)$, we can define the so-called
``continuous interpolation spaces" of the couple $(E_1, E_0)$.
This way of defining $E_{\theta}$ is described in \cite{DP1}. For
any $v \in E_1(J)$ we have $|\dot{v}|_{E_0}\leq C t^{1-\theta}$
for some finite constant $C>0,$ so that $v$ extends to a
continuous function from $[0,1]$ to $E_0$. Indeed $v(0)$ is well
defined by \begin{eqnarray*}
v(0)=v(1)-\int_0^1\dot{v}(t)dt.\end{eqnarray*} We can define the
Banach space \begin{eqnarray*} E_{\theta}=\left\{u(0):u \in
E_1(J)\right\},~~~~|x|_{E_{\theta}}=\inf \left\{|u|_{E_1(J)}:
x=u(0), u\in E_1(J)\right\}.\end{eqnarray*}

Let $E_{\theta}, E_0(J), E_1(J)$ as above. Suppose that $A \in$
$\cal L$$(E_1, E_0)$, we can consider the operator $\hat{A}:
E_1(J)\rightarrow E_0(J)\oplus E_{\theta}$ defined by
$\hat{A}u=(\dot{u}-Au, u(0))$. $\hat{A}$ is obviously bounded. We
define $M_\theta(E_1, E_0)$ as follows: $M_\theta(E_1, E_0)=\{A
\in $ $\cal H$ $(E_1, E_0): \hat{A}$ is an isomorphism between
$E_1(J)$ and $E_0(J)\oplus E_{\theta}\}.$ If $A \in M_\theta(E_1,
E_0)$ then $(E_0, E_1)$ is called \textbf{a pair of maximal
regularity} for A. In other words, $M_\theta(E_1, E_0)$ contains
those generators of analytic semigroups $A$ for which
\begin{eqnarray*}\dot{u}+Au(t)&=&f(t),\\
u(0)&=&u_0,\end{eqnarray*} has a unique solution $u \in E_1(J)$
for any $f \in E_0(J)$, and $u_0 \in E_{\theta}$. In \cite{An}, S.
Angenent has shown that many interesting operators belong to the
class $M_\theta(E_1, E_0)$. In particular the operator
$A(\varphi)$ in (\ref{3-1}) belongs to the class $M_\theta(E_1,
E_0)$ for appropriate Banach space $E_0, E_1$ (see next section).\\

Now we are in the position to state a main Theorem in \cite{CS1}.
\begin{theo}[Philippe Cl\'ement, Gieri Simonett \cite{CS1}]
Let $\theta \in (0, 1)$ be fixed and let $E_\theta=(E_0,
E_1)_{\theta}$ be a continuous interpolation space. Assume
$V_\theta \subset E_\theta$ is open,
\begin{eqnarray*}(A, f) \in
C^{1-}\left(V_\theta, M_\theta(E_1, E_0)\times
E_0\right).\end{eqnarray*} For every $x_0 \in V_\theta$, there
exists some positive constant $\tau=\tau(x_0),
\epsilon=\epsilon(x_0)$ and $c=c(x_0)$ such that the following
evolution equation
\begin{eqnarray*}\dot u(t)+A(u(t))u(t)&=&f(u(t)), t \in
\widetilde{J},\\
u(0)&=&x,\end{eqnarray*} has a unique solution \begin{eqnarray*}
u(.,x)\in C^1_{1-\theta}\left([0, \tau], E_0\right)\cap
C_{1-\theta}\left([0, \tau], E_1\right)\end{eqnarray*} in $[0,
\tau]$ for any initial value $x\in \bar{B}_{E_\theta}(x_0,
\epsilon)$, where $\bar{B}_{E_\theta}(x_0, \epsilon)$ denotes the
ball in Banach space $E_\theta$ with radius $\epsilon$ centered at
$x_0$. In particular, we have \begin{eqnarray*} |u(\cdot,
x_0)|_{C([0, \tau], E_\theta)}\leq c,~\mbox{and}~|u(\cdot,
x_0)|_{E_1([0,\tau])}\leq c.\end{eqnarray*} Moreover, we
have\begin{eqnarray*}\left|u(\cdot,x)-u(\cdot,y)\right|_{E_1([0,
\tau])}&\leq& c|x-y|_{E_\theta},\end{eqnarray*} and
\begin{eqnarray*} \left|u(\cdot,x)-u(\cdot,y)\right|_{C([0, \tau],
E_\theta)}\leq c|x-y|_{E_\theta},~~x, y \in
\bar{B}_{E_\theta}(x_0, \epsilon) .\end{eqnarray*}
\end{theo}

\subsection{Short time existence}
~~~~We now exhibit appropriate Banach space $E_0, E_1$ to use
Theorem 3.1 to prove the short time existence of the Calabi flow.
The spaces we use are certain {\it little H\"older space}.  Recall
that if $k \in \mathbb{N}, \alpha \in (0, 1)$, the H\"older space
$C^{k, \alpha}$ is the Banach space of all $C^k$ functions $f:
\mathbb{R}^n\rightarrow \mathbb{R}$ which has finite H\"older
norm. The subspace of $C^{\infty}$ function in $C^{k, \alpha}$ is
not dense under the H\"older norm. One defines the little H\"older
space $c^{k, \alpha}$ to be the closure of smooth functions in the
usual H\"older space $C^{k,\alpha}$. And one can verify $c^{k,
\alpha}$ is a Banach space and that $c^{l,\beta}\hookrightarrow
c^{k, \alpha}$ is a continuous and dense imbedding for $k\leq l$
and $0<\alpha<\beta<1$. These definitions can be extended to
functions on a smooth manifold $M$ naturally \cite{BJ}. For our
purpose, the key fact \cite{Tr} about the continuous interpolation
spaces is that for $k\leq l$ and $0<\alpha<\beta<1$, and
$0<\theta<1$, there is a Banach space isomorphism
\begin{equation} \label{3-2}(c^{k, \alpha}, c^{l, \beta})_{\theta}\cong
c^{\theta l+(1-\theta)k+\theta
\beta+(1-\theta)\alpha},\end{equation} provided that the exponent
$\theta l+(1-\theta)k+\theta \beta+(1-\theta)\alpha$ is not an
integer.

For any function $\varphi\in \cal H$$_\omega$, the operator
$A(\varphi)$ is a strict elliptic operator. Under these
conditions, $A$ generates an analytic semigroup \cite{Lu} in
$c^{0,\alpha}(M)$ with domain $c^{4,\alpha}(M).$  By Theorem 3.1,
the short time existence holds for any initial potential $\varphi
\in c^{3, \alpha}(M)$ satisfying $\omega_{\varphi}>0$. We have

\begin{theo}
If $\omega_{\varphi_0}=\omega+\i\p\bar{\p}\varphi_0$ satisfies
$|\varphi_0|_{c^{3, \alpha}(M)}<K$, where $K$ is some constant,
and $\lambda \omega< \omega_0=\omega_{\varphi_0}< \Lambda \omega$,
where $\lambda, \Lambda$ are two positive constants, then the
Calabi flow initiating with $\varphi_0$ admits a unique solution
\begin{eqnarray*}\varphi(t) \in C([0,T], c^{3, \alpha}(M))\cap C((0, T], c^{4,
\alpha}(M))\end{eqnarray*} for some small $T=T(\lambda, \Lambda,
K, g)$. More specifically, for any $t \in (0, T]$, there is a
constant $C=C(\lambda, \Lambda, K, g)$ such that
\[
t^{1/4}\left(|\dot{\varphi}(t)|_{c^{0,\alpha}(M)}+|\varphi(t)|_{c^{4,
\alpha}(M)}\right)\leq C,~~|\varphi(t)|_{c^{3, \alpha}(M)}\leq C.
\]
In particular, when $g$ is CscK, and $|\varphi_0|_{c^{3,
\alpha}(M)}$ is small enough, there is some uniform constant $C_1$
depending only on $g$, such that \begin{eqnarray*}
|\varphi(t)|_{c^{3, \alpha}(M)}\leq C_1|\varphi_0|_{c^{3,
\alpha}(M)},\end{eqnarray*} and \begin{eqnarray*}
t^{1/4}\left(|\dot{\varphi}(t)|_{c^{0,\alpha}(M)}+|\varphi(t)|_{c^{4,
\alpha}(M)}\right)&\leq& C_1|\varphi_0|_{c^{3,\alpha}(M)}.
\end{eqnarray*}
\end{theo}


\begin{proof}
Set $E_0=c^{0,\alpha}(M), E_1=c^{4,\alpha}(M)$ for any fixed
$\alpha \in (0,1)$. Choose $\theta=3/4$ so that $E_\theta=(E_0,
E_1)_\theta=c^{3, \alpha}(M)$. For any $\varphi_0 \in V_\theta$,
where \begin{eqnarray*} V_\theta:=\left\{u|u\in c^{3,\alpha}(M),
\lambda g<g(u)<\Lambda g,
|u|_{c^{3,\alpha}(M)}<K\right\}\end{eqnarray*} is open in
$c^{3,\alpha}(M)$, $A(\varphi_0)=A$ generates a strongly
continuous analytic semigroup $\exp(-tA)$ in $c^{0, \alpha}(M)$
with domain $c^{4, \alpha}(M)$ and $A \in M_\theta(E_1, E_0)$.
This fact follows from the construction in \cite{An}. Because for
any $0<\alpha<1, A \in \cal H$$(E_1, E_0).$ We can take $F_0=c^{0,
\gamma}(M), F_1=c^{4, \gamma}(M)$ for some $0<\gamma<\alpha<1,$
and $A \in \cal H$$(F_1, F_0).$ Then define $F_2$ to be the domain
of $A^2$, i.e. \begin{eqnarray*} F_2=\{x \in F_1: Ax \in
F_1\},\end{eqnarray*} and let $F_2$ have the graph topology,
$|x|_{F_2}=|x|_{F_1}+|Ax|_{F_1}$. In our case $F_2=c^{8,
\gamma}(M)$. Pick up $\delta=(\alpha-\gamma)/4$ such
that\begin{eqnarray*} E_0=F_{\delta}, ~~~E_1=F_{1+\delta}=(F_2,
F_1)_{\delta}.\end{eqnarray*} Then by Theorem 2.14 in \cite{An} we
know actually $A \in
M_{\theta}(E_1, E_0)$.\\

To use Theorem 3.1, we need to show $(A, f) \in C^{1-}(E_\theta,
M_\theta(E_1, E_0)\times E_0)$ for $\theta=3/4$, and
$E_0=c^{0,\alpha}(M), E_1=c^{4,\alpha}(M).$ It suffices to show
that for $u, v \in V_\theta, w\in E_1= c^{4,\alpha}(M),$
\begin{eqnarray*}
|f(u)-f(v)|_{c^{0,\alpha}(M)}&&\leq C_1|u-v|_{c^{3,\alpha}(M)}, \\
|A(u)w-A(v)w|_{E_0}&&\leq
C_2|u-v|_{c^{3,\alpha}(M)}|w|_{E_1},\end{eqnarray*} where $C_1,
C_2$ are two constants. The first inequality is obvious since $f$
is actually $C^{\infty}$ in its arguments. The second one follows
from \begin{eqnarray*}
\left|A(u)w-A(v)w\right|_{E_0}&=&\left|[g^{i\bar j}(u)g^{k\bar
l}(u)-g^{i\bar j}(v)g^{k\bar l}(v)] \nabla_{i}\nabla_{\bar
j}\nabla_{k}\nabla_{\bar l}
         w\right|_{E_0}\\
&\leq&C|u-v|_{c^{2,\alpha}(M)}|w|_{E_1}\\
&\leq&C|u-v|_{c^{3,\alpha}(M)}|w|_{E_1}.\end{eqnarray*} Theorem
3.2 follows from Theorem 3.1 directly.
\end{proof}
\begin{theo}(Smoothing property)
Following Theorem $3.2$, actually $\varphi(t)\in C\left([0, T],
c^{3, \alpha}(M)\right)\cap C\left((0, T], C^{\infty}(M)\right)$.
\end{theo}
\begin{proof}Consider the following
parabolic equation \begin{eqnarray*}\frac{\p \varphi}{\p
t}&=&-A(\varphi)\varphi+f(\varphi)\\
\varphi(0)&=&\varphi_0.\end{eqnarray*} By the maximal regularity,
we know indeed $\varphi(t) \in c^{4, \alpha}(M)$ when $t>0$. We
can consider \begin{eqnarray*}\frac{\p \varphi}{\p
t}&=&-A(\varphi)\varphi+f(\varphi),\end{eqnarray*} with initial
value $\varphi(\delta) \in c^{4, \alpha}(M)$ for any $\delta>0$
small. Taking $E_0=c^{1, \alpha}(M), E_1=c^{5, \alpha}(M)$, by
(\ref{3-2}) we have $E_{\theta}=(E_0, E_1)_{\theta}=c^{4,
\alpha}(M)$ where $\theta=3/4.$ Theorem 3.1 applies for $(E_0,
E_1)$, a pair of maximal regularity for $A(\varphi(\delta))$. So
we can get $\varphi(t) \in c^{5, \alpha}(M)$ when $t \in (\delta,
T].$ Similarly we can show $\varphi(t) \in c^{k, \alpha}(M)$ for
any $k \in \mathbb{N}$ when $t>0$.
\end{proof}

\begin{rem}One can show $\varphi(t) \in C^{\infty}\left( (0,
T], C^{\infty}(M)\right).$ To show $\varphi(t)$ is smooth in time
argument, we need to consider some Banach space with time norm
\cite{Lu}. The idea is similar to Theorem 3.3, but we shall not
use this fact.
\end{rem}

\section{Stability theorem}
~~~~When $g$ is a CscK metric, it is a fixed point under the
Calabi flow. Since $\cC a(\omega_\varphi)$ is weakly convex at
CscK metric $g$, intuitively it is a stable fixed point under the
Calabi flow. We will give an affirmative answer in this section.
More precisely,
\begin{theo} (Stability theorem)
Suppose $g$ is a CscK metric in $[\omega]$ on $M$. Consider the
Calabi flow
\begin{equation}
  \frac{\p \varphi}{\p t}=R_\varphi-\underline{R},
  \qquad
  \varphi(0)=\varphi_0.
\end{equation}
When $|\varphi_0|_{c^{3, \alpha}(M)}<\epsilon$, where $\epsilon>0$
is small enough, then the flow exists for all time and
$g(t)=g_{\varphi(t)}$ converges to a CscK $g_\infty$ metric in the
same class $[\o]$ in $C^{\infty}$ sense. Moreover, the Calabi
energy decays exponentially fast and the convergence is
exponentially fast.
\end{theo}
\begin{proof}
Suppose $g$ is a CscK metric, we derive {\it a priori} estimates
for the scalar curvature equation. Consider the following elliptic
equation: \begin{eqnarray*}
R_\varphi-\underline{R}=f,\end{eqnarray*} where
$\omega_\varphi=\omega+\i\p \bar{\p}\varphi$ is ``$C^{\alpha}$
equivalent" to $\omega$ for some $\alpha \in (0, 1)$. By
``$C^{\alpha}$ equivalent" we mean in a holomorphic local
coordinates,  $\lambda \omega\leq \omega_\varphi=\omega+\i\p
\bar{\p}\varphi \leq \Lambda \omega$ for some constants $0<\lambda
\leq \Lambda$, and $|\varphi|_{C^{2, \alpha}(M, g)}$ is uniformly
bounded. Rewrite the equation as
\begin{eqnarray*} -g^{i\bar{j}}_{\varphi}\p_i\p_{\bar{j}}\log
\det{(g_{k\bar{l}}+\p_k\p_{\bar{l}}\varphi)}+g^{i\bar{j}}\p_i\p_{\bar{j}}\log
\det{(g_{k\bar{l}})}=f.\end{eqnarray*} It follows that
\begin{eqnarray}\label{4-2}
-g^{i\bar{j}}_{\varphi}\p_i\p_{\bar{j}}\log\frac{\det{(g_{k\bar{l}}+\p_k\p_{\bar{l}}\varphi)}}
{\det{(g_{k\bar{l}}})}=f+\left(g^{i\bar{j}}_{\varphi}-g^{i\bar{j}}\right)\p_i\p_{\bar{j}}\log\det{(g_{k\bar{l}})}
.\end{eqnarray} Denote \begin{eqnarray*}
u&=&\log\frac{\det{(g_{k\bar{l}}+\p_k\p_{\bar{l}}\varphi)}}{\det{(g_{k\bar{l}}})},\\
h&=&f+\left(g^{i\bar{j}}_{\varphi}-g^{i\bar{j}}\right)\p_i\p_{\bar{j}}\log\det{(g_{k\bar{l}})}
.\end{eqnarray*} One can rewrite (\ref{4-2}) as
\begin{equation}\label{4-3}
-\triangle_\varphi u=h.
\end{equation}
By the standard $L^p$ theory ($p$ big enough), one can get that
\begin{equation}\label{4-4} |u-\underline{u}|_{W^{2, p}}\leq
C(|h|_{L^p}+|u-\underline{u}|_{L^p}),
\end{equation}
where $\underline{u}$ is the average of $u$, the constant $C$
depends on the $C^{\alpha}$ H\"older norm of $\o_{\varphi}$
(actually $C^{0}$ norm is sufficient for $L^{p}$ theory).
 Through Moser's iteration the estimate
\[|u-\underline{u}|_{L^\infty}\leq C(|h|_{L^p}+|u-\underline{u}|_{L^2})
\] holds for (\ref{4-3}) when $g_\varphi$ is equivalent to $g$.
And by (\ref{4-3}), one can get
\[
-(u-\underline{u})\triangle_\varphi(u-\underline{u})=h(u-\underline
{u}).
\]
It implies that
\[
\int_M|\nabla
u|^2dg_\varphi=\int_Mh(u-\underline{u})dg_\varphi\leq
\left(\int_Mh^2dg_\varphi\right)^{1/2}\left(\int_M(u-\underline{u})^2dg_\varphi\right)^{1/2}.
\]
The Poincar\'e inequality gives that
\[
\int_M|u-\underline{u}|^2dg_\varphi\leq C\int_M|\nabla
u|^2dg_\varphi.
\]
And so \[|u-\underline{u}|_{L^2}\leq C|h|_{L^2}.
\]
It implies that
\[
|u-\underline{u}|_{L^\infty}\leq C(|h|_{L^p}+|h|_{L^2})\leq
C|h|_{L^p}.
\]
We need only the weaker version
\begin{equation}\label{4-5}
|u-\underline{u}|_{L^p}\leq C|h|_{L^p}.
\end{equation}
It follows from (\ref{4-4}) and (\ref{4-5}) that
\begin{equation}\label{4-6}
|u-\underline{u}|_{W^{2, p}}\leq C|h|_{L^p}.
\end{equation}

We can rewrite (\ref{4-6}) as \begin{eqnarray}\label{4-7}
|u-\underline{u}|_{W^{2,p}}&\leq&
C\left|f+\left(g^{i\bar{j}}_{\varphi}-g^{i\bar{j}}\right)\p_i\p_{\bar{j}}\log\det{(g_{k\bar{l}})}\right|_{L^{p}}\nonumber\\
&\leq&C\left(|f|_{L^p}+|\varphi|_{W^{2, p}}\right).\end{eqnarray}
 Consider the Monge-Amp\'ere equation, \begin{eqnarray*}
\log\frac{\det{(g_{i\bar{j}}+\p_i\p_{\bar{j}}\varphi)}}{\det{(g_{i\bar{j}}})}=u.
 \end{eqnarray*}
We can rewrite it as
\[
\frac{\det{(g_{i\bar{j}}+\p_i\p_{\bar{j}}\varphi)}}{\det{(g_{i\bar{j}}})}=\exp{(\underline{u})}\exp{(u-\underline{u})}
\]
 Since $g_\varphi$ is ``$C^{\alpha}$ equivalent" to $g$, we know $\exp{(\underline{h})}$ is a uniformly bounded constant and so when $p$ is big ($p>2n$)
 \begin{equation}\label{4-8} \left|\varphi\right|_{C^{3, \alpha}}\leq C |u-\underline{u}|_{C^{1, \alpha}}\leq C
 |u-\underline{u}|_{W^{2,p}}\end{equation} for some $\alpha \in (0,1).$
Combine (\ref{4-7}) and (\ref{4-8}), we have \begin{eqnarray*}
|\varphi|_{C^{3, \alpha}}\leq
C\left(|f|_{L^p}+|\varphi|_{W^{2,p}}\right).\end{eqnarray*} By the
interpolation inequality, actually we have
\begin{equation}\label{4-9} |\varphi|_{C^{3,\alpha}}\leq
C\left(\left|R_\varphi-\underline{R}\right|_{L^p}+|\varphi|_{L^1}\right).\end{equation}
Since the metrics are equivalent, it does not matter which volume
form we use as measure. Also $C$ denotes a generic constant which
could differ line by line.

To derive the long time existence, it is sufficient to control the
$C^{3, \alpha}$ norm of the potential by the short time existence
derived in Section 3. By (\ref{4-9}), essentially  we just  need
to control the decay of the Calabi energy. Indeed, the Calabi
energy decays exponentially fast around a CscK metric. Let us
recall
\begin{eqnarray*}
\frac{d}{dt}\int_M(R(t)-\underline{R})^2dg(t)=-2\int_M\left(D_tR(t),
R(t)\right) dg(t).\end{eqnarray*} By the definition of
eigenvalues, the first eigenvalue of $D_t$ is given by
\begin{eqnarray*} \lambda(t)=\inf_{f \in C_0^{\infty}(M)}\frac{\int_M
f,_{\alpha\beta}f^{,\alpha\beta}dg(t)}{\int_Mf^2dg(t)}.\end{eqnarray*}
The definition involves only the metric and its first derivative.
If $g(t)$ is close to $g$ in $C^{3, \alpha}-$topology (by that we
mean $\varphi(t)$ is small in $C^{3, \alpha}$,
$|\varphi(t)|_{C^{3,\alpha}}<\epsilon$), then the eigenvalues of
$D_t$ is close to the eigenvalues of $D$, and the corresponding
eigenspace is close to the eigenspace of $g$. The eigenvalues of
$D$ are nonnegative.\\

\noindent{\bf Case 1.} If the first eigenvalue of $D$ is strictly
bigger than zero, then for any $D_t$, the first eigenvalue is
bounded away from zero. It implies
\begin{eqnarray*}\int_M\left(D_tR(t), R(t)\right) dg(t)\geq \delta
\int _M
\left(R(t)-\underline{R}\right)^2dg(t),\end{eqnarray*} for some constant $\delta>0.$ \\

\noindent{\bf Case 2.} If the first eigenvalue of $D$ is zero, by
using the Futaki character, we can still prove the above estimate
\cite{Chen}. Notice that when the K\"ahler class admits a CscK
metric, the Futaki character vanishes. It implies
\begin{eqnarray*}\int_MXF(t)dg(t)=0,\end{eqnarray*} for any $X \in \h(M)$. Rewrite the
Futaki character as the following,
\begin{eqnarray*}\int_M\theta_X(R(t)-\underline{R})dg(t)=0,\end{eqnarray*} where $\theta_X$ is the
potential of the vector field $X$ with respect to the metric
$g(t).$  The potential function $\theta_X$ of a holomorphic vector
field $X$ is defined (up to addition of some constant) as the
following: \begin{eqnarray*}\int_M
X(f)dg=-\int_M\theta_X\triangle_gfdg.\end{eqnarray*} The
eigenvalue of $D_t$ converges to the eigenvalue of $D$, and the
eigenspace $\Upsilon_t$ of the first eigenvalue of $D_t$ converges
to the first eigenspace $\Upsilon$ of $D$ when $\epsilon$ goes to
zero. And if the first eigenvalue of $D$ is zero, then the first
eigenspace $\Upsilon$ consists of functions which are the
potential of some holomorphic vector fields (if $f \in \Upsilon$,
$\nabla f$ is the real part of a holomorphic vector field). If we
decompose $R-\underline{R}$ into two components, then
\begin{eqnarray*} R-\underline{R}=\rho+\rho^\perp,\end{eqnarray*} where
$\rho \in \Upsilon_t$ and $\rho^\perp \perp \Upsilon_t.$ The
vanishing of the Futaki character implies that
\begin{eqnarray*}\parallel\rho\parallel_{L^2(g(t))}\leq
c(\epsilon)\parallel
R(t)-\underline{R}\parallel_{L^2(g(t))},\end{eqnarray*} where
$c(\epsilon)\rightarrow 0$ as $\epsilon \rightarrow 0.$ If
$\epsilon$ is small enough, we can get \begin{eqnarray*}
\int_M\left(D_tR(t), R(t)\right) dg(t)&\geq&
\lambda_1(t)\parallel\rho\parallel^2_{L^2(g(t))}+\lambda_2(t)\parallel
\rho^\perp\parallel^2_{L^2(g(t))}\\
&\geq& \frac{\lambda_2}{2}\parallel
R(t)-\underline{R}\parallel^2_{L^2(g(t))},\end{eqnarray*} where
$\lambda_i(t)$ is the i-th eigenvalue of $D(t)$, and $\lambda_i$
is the i-th eigenvalue
of $D$. By assumption we know $\lambda_1=0, \lambda_2>0.$ \\

In both cases the Calabi energy exponentially decays:
\begin{eqnarray*} \int_M\left(R(t)-\underline{R}\right)^2dg(t)\leq \exp(-\delta t)
\int_M\left(R(0)-\underline{R}\right)^2dg(0),\end{eqnarray*} where
$\delta$ is a positive
constant bounded away from zero.\\

Now we are in the position to derive the long time existence and
the convergence. By the smoothing property of the Calabi flow
(Theorem 3.3), we can actually assume $\varphi_0 \in C^{\infty}$
and $|\varphi_0|_{C^{4, \alpha}}< \epsilon.$ Define
\begin{eqnarray*} \cG=\left\{\varphi:\quad \varphi \in
C^{\infty},\quad \lambda g<g_\varphi<\Lambda g,\quad
|\varphi|_{C^{3, \alpha}}< K\right\}.\end{eqnarray*} If
$\epsilon>0$ is small, $\varphi_0 \in \cG$, by Theorem 3.2, we
know the flow exists at least in $[0, T]$, where $T$ depends on
$\cG.$ Moreover, $|\varphi(t)|_{C^{4, \alpha}}\leq C
|\varphi_0|_{C^{4, \alpha}}$ for any $t \in [0, T]$, where $C$ is
a constant depending only on $\cG$. It implies in $[0,T]$, the
Calabi energy decays exponentially fast because $g(t)$ is close to
$g$ in $C^{4, \alpha}$.  And we know \begin{eqnarray*}
\varphi(t)=\varphi(0)+\int_0^t(R(s)-\underline{R})ds,
\end{eqnarray*} so we can get
\begin{eqnarray*}\int_M|\varphi(t)|\omega^n\leq
\int_M|\varphi(0)|\omega^n+\int_M\left|\int_0^t(R(s)-\underline{R})ds\right|\omega^n.\end{eqnarray*}
Note that in $[0, T]$, all the metrics are close to CscK metric
$g$ in $C^{4, \alpha}$, in particular all the metrics are
``$C^{\alpha}$ equivalent" to $g$. By Schwartz's inequality, we
can get
\begin{eqnarray*}\int_M\left|\int_0^t(R(s)-\underline{R})ds\right|\omega^n&\leq&
\int_0^t\int_M\left|R(s)-\underline{R}\right|\omega^nds\\
&\leq&C\int_0^t\int_M\left|R(s)-\underline{R}\right|\omega^n(s)ds\\
&\leq&C\int_0^t\left(\int_M\left(R(s)-\underline{R}\right)^2\omega^n(s)\right)^{1/2}\left(\int_M\omega^n(s)\right)^{1/2}ds\\
&\leq&C\left(\int_M\left(R(0)-\underline{R}\right)^2\omega^n(0)\right)^{1/2}\int_0^t\exp(-cs/2)ds\\
&\leq&C\left|R(0)-\underline{R}\right|_{L^2}.\end{eqnarray*} Since
$g(t)$ is close to $g$ in $C^{4, \alpha}$, we can assume
$|R(t)-\underline{R}|\leq 1$, then
\begin{eqnarray*}\left|R_\varphi-\underline{R}\right|_{L^p}\leq
C\left|R_\varphi-\underline{R}\right|^{2/p}_{L^2}\leq
C|R(0)-\underline{R}|^{2/p}_{L^2}.\end{eqnarray*} By (\ref{4-9}),
we know $\forall t \in [0, T]$\begin{eqnarray*}
|\varphi(t)|_{C^{3, \alpha}}&\leq&
C\left(\left|R_\varphi-\underline{R}\right|_{L^p}+|\varphi(t)|_{L^1}\right)\\
&\leq&C\left(\left|R(0)-\underline{R}\right|^{2/p}_{L^2}+|\varphi(0)|_{L^1}+|R(0)-\underline{R}|_{L^2}\right),
\end{eqnarray*} where $C$ is a uniform constant independent of time $t$. When
$\epsilon
>0$ is small enough, \begin{eqnarray*}|\varphi(T)|_{C^{3, \alpha}}\leq C
c(\epsilon),\end{eqnarray*} where $c(\epsilon)\rightarrow 0$ when
$\epsilon \rightarrow 0$.  Then we can extend the flow to $[T,
2T]$ starting from $\varphi(T) \in \cG$. And in $[0, 2T]$, all the
metrics are close to CscK metric $g$. By (\ref{4-9}), we still
have $\forall t \in [0, 2T]$,
\begin{eqnarray*}|\varphi(t)|_{C^{3,\alpha}}&\leq&
C\left(\left|R_\varphi-\underline{R}\right|_{L^p}+|\varphi(t)|_{L^1}\right)\\
&\leq&C\left(|\varphi(0)|_{L^1}+|R(0)-\underline{R}|^{2/p}_{L^2}+|R(0)-\underline{R}|_{L^2}\right)\\
&\leq& Cc(\epsilon).\end{eqnarray*} Similarly we can extend the
flow to $[2T, 3T]$ starting from $\varphi(2T)$ with the uniform
bounded $C^{3, \alpha}$ norm, $\forall t \in [0, 3T]$
\begin{eqnarray*}
|\varphi(t)|_{C^{3,\alpha}}&\leq&C\left(|\varphi(0)|_{L^1}+|R(0)-\underline{R}|^{2/p}_{L^2}+|R(0)-\underline{R}|_{L^2}\right)\\
&\leq& Cc(\epsilon).\end{eqnarray*} By iteration the flow exists
for all time and actually, $\varphi(t)$ is uniformly small in
$C^{3, \alpha}$. By smoothing property of the Calabi flow, indeed
$g(t)$ is uniformly close to $g$ in $C^{\infty}$ and the Calabi
energy decays exponentially fast. In particular, $\varphi(t)$ is
uniformly bounded in $C^{\infty}$ norm. And also we know that for
any $t, \bar t\in (0, \infty)$
\begin{eqnarray}\label{4-10}|\vphi(t)-\vphi(\bar t)|_{L^1}&=&\int_M|\vphi(t)-\vphi(\bar t)|dg
\nonumber\\&=&\int_M\int_{t}^{\bar t}\frac{\p}{\p s}\vphi(s)dsdg\nonumber\\
&\leq&C\int_{t}^{\bar t}\int_M|R(\vphi(s))-\underline{R}|dgds\nonumber\\
&\leq&C\int_t^{\bar t}\left(|R(\vphi(s))-\underline{R}|_{L^2}\right)ds\nonumber\\
&\leq&C|R(0)-\underline{R}|_{L^2}\int_t^{\bar t}\exp(-\delta
s/2)ds.
\end{eqnarray}
So there exists a limit potential $\varphi(\infty)\in C^{\infty}$,
such that $\varphi(t)$ converges to $\varphi(\infty)$ by sequence
in $C^{\infty}$ and in $L^1$ in the flow sense. In particular,
$\varphi(\infty)$ defines a CscK metric $g_\infty$ in the same
class $[\o].$  Also we can get  {\it a priori} estimates for
\[R(\vphi(t))-R(\vphi(\infty))=f.\]
By the exact same argument, one can get
\begin{equation}\label{4-11}|\vphi(t)-\vphi(\infty)|_{C^{3, \alpha}(M)}\leq
C\left(|f|_{L^p}+|\vphi(t)-\vphi(\infty)|_{L^1}\right).\end{equation}
By (\ref{4-10}) and (\ref{4-11}),  the flow converges not only by
sequence but exponentially fast along the flow.
\end{proof}
\begin{rem}The limit CscK metric $g_\infty$ could be different
with  $g$ in holomorphic coordinates. It means in general
$\varphi(\infty)\neq 0.$ To see this, we can pick up
$\varphi(0)\neq 0$ such that $g_0$ is a CscK metric and is
sufficient close to $g$ in holomorphic coordinates (but $g_0\neq
g$). But along the flow, $\varphi(t)\equiv\varphi(0)\neq 0.$  For
example, on Riemann sphere $S^2$, $(S^2,
\lambda^2|dz|^2/(1+\lambda^2|z^2|)^2)$ is a smooth continuous
family of constant Gauss curvature metrics for $\lambda\in (0,
\infty)$. From technical point of view, it explains also that the
$L^1$ norm of the potential in (\ref{4-9}) is necessary.
\end{rem}
\section{Long time existence}
~~~~Essentially, to understand the whole picture of the Calabi
flow, we need to understand some compactness behavior of
$\omega_\varphi \in \cH_\omega$ with bounded Calabi energy since
the Calabi flow is the gradient flow which decreases the Calabi
energy. In this section, we derive a compactness theorem in
$\cH_\omega$. As a consequence, the long time existence holds
under the assumption of uniformly bounded Ricci curvature.
\begin{theo}(Compactness theorem) All metrics
$\omega_\varphi=\omega+\sqrt{-1}\p \bar{\p }\varphi$ in the space
of K\"ahler metrics $\cH_\o$ with both the potential $\varphi$ and
the Ricci curvature $Ric_\varphi$ uniformly bounded are equivalent
and compact in $C^{1, \alpha}-$topology for any $\alpha \in (0,1)$
(It is equivalent to say $\varphi$ is uniformly bounded in $C^{3,
\alpha}$ for any $\alpha \in (0,1)$).
\end{theo}


\begin{proof}  Denote \begin{eqnarray*} F=\log
\left(\frac{\omega^n_\varphi}{\omega^n}\right).\end{eqnarray*}  We
know
\begin{eqnarray*}
\triangle F&=&g^{i\bar{j}}\p_i\p_{\bar{j}}\log\left(\frac{\omega^n_\varphi}{\omega^n}\right)\\
&=&-g^{i\bar{j}}R_{i\bar{j}}(\varphi)+R ,\end{eqnarray*} where
$\triangle$ and scalar curvature $R$ are both with respect to the
background metric $\omega$, and $R_{i\bar{j}}(\varphi)$ is the
Ricci curvature of $\omega_\varphi$. Since $|Ric_\varphi|$ is
bounded, then \begin{eqnarray*}\left|\triangle F-R\right|\leq
C(n+\triangle \varphi).\end{eqnarray*} In other words,
\begin{equation}\label{5-1} \triangle (F-C\varphi)\leq C_1,\end{equation} and
\begin{equation}\label{5-2}\triangle\left( F+C\varphi\right)\geq-C_2\end{equation} for
some positive constant $C_1, C_2.$ By (\ref{5-1}) (\ref{5-2}), we
can deduce the $C^0$ bound of $F$. For simplicity, set the volume
of $M$ to be $1$,
\begin{eqnarray*}\int_M\omega^n=1.\end{eqnarray*} Let us recall the Green's formula on a
compact manifold,\begin{eqnarray*}
\left(F+C\varphi\right)(p)=-\int_MG(p,q)\left\{\triangle\left(F+C\varphi\right)(q)\right\}\omega^n(q)
+\int_M\left(F+C\varphi\right)\omega^n ,\end{eqnarray*} where
$G(p,q)\geq 0$ is Green's function of $\omega$. The first term of
right hand side is bounded from above since
$\triangle\left(F+C\varphi\right)$ is bounded from below. To bound
the second term, we have
\begin{eqnarray*}\exp\left(\int_MF\omega^n\right)\leq\int_M\exp(F)\omega^n=1\end{eqnarray*}
provided that \begin{eqnarray*}\int_M\omega^n=1.\end{eqnarray*}
Because $\varphi$ is uniformly bounded, we get that $F$ is bounded
from above. Before we derive the lower bound of $F$, we try to
show that $\triangle \varphi$ is uniformly bounded following Yau's
celebrated work on Calabi conjecture. Consider the following
Monge-Ampere equation
\begin{eqnarray*}\frac{\det(g_{i\bar{j}}+\varphi,_{i\bar{j}})}{\det(g_{i\bar{j}})}=\exp(F).\end{eqnarray*}
Just following Yau's work \cite{Yau}, for some big constant $C_3$,
we get at some point $p$ (the maximum of
$\exp(-C_3\varphi)(n+\triangle \varphi)$),\begin{eqnarray*}
\triangle_{\varphi}\left\{\exp(-C_3\varphi)(n+\triangle
\varphi)\right\}(p)\leq0.\end{eqnarray*} Follow Yau's calculation,
at the point $p$,\begin{eqnarray}\label{5-3}0&\geq& \triangle
F-n^2\inf_{i\neq
l}R_{i\bar{i}l\bar{l}}-C_3n(n+\triangle \varphi)\nonumber\\
&&\quad+\left(C_3+\inf_{i\neq
l}R_{i\bar{i}l\bar{l}}\right)\exp\left\{\frac{-F}{n-1}\right\}(n+\triangle
\varphi)^{n/(n-1)} .\end{eqnarray} By (\ref{5-2}) (\ref{5-3}), at
the point $p$ we have, \begin{eqnarray}\label{5-4}0&\geq&
-C\triangle \varphi-C_2-n^2\inf_{i\neq
l}R_{i\bar{i}l\bar{l}}-C_3n(n+\triangle \varphi)\nonumber\\
&&\quad+\left(C_3+\inf_{i\neq
l}R_{i\bar{i}l\bar{l}}\right)\exp\left\{\frac{-F}{n-1}\right\}(n+\triangle
\varphi)^{n/(n-1)}.\end{eqnarray} (\ref{5-4}) implies that
$(n+\triangle\varphi)(p)$ has an upper bound $C_0$ depending only
on $\sup_MF$ and $M$. It follows that
\begin{eqnarray*} \exp(-C_3\varphi)(n+\triangle \varphi) &\leq&
\exp(-C_3\varphi)(n+\triangle \varphi)(p)\\
&\leq& C_0\exp(-C_3\varphi(p)).\end{eqnarray*} Because $\varphi$
is uniformly bounded, we obtain \begin{equation}\label{5-5}
0<n+\triangle \varphi \leq C(C_3, \sup_MF, M).\end{equation}So
$\triangle \varphi$ is uniformly bounded. It follows from
(\ref{5-1}) (\ref{5-2}) (\ref{5-5})
that  $\triangle F$ is uniformly bounded .\\

The lower bound of $F$ follows from the following result, which
is essentially the proposition 3.14 in \cite{CT}:
\begin{prop}\cite{CT} If Ricci curvature is bounded from below and the
K\"ahler potential is uniformly bounded.   Then, there is a
uniform constant $C$ such that:\begin{eqnarray*}
\inf_MF\geq-4C\exp\left(2+2\int_M\log\frac{\omega^n_{\varphi}}{\omega^n}\omega^n_{\varphi}\right).
\end{eqnarray*}
\end{prop}

\begin{proof} Choose a constant $c$ such that
\begin{equation*}\int_M\log\frac{\omega^n_{\varphi}}{\omega^n}\omega^n_{\varphi}\leq
c.\end{equation*} And remember we set
\[
\int_M\omega^n=1.
\]
Put $\epsilon$ to be $\exp(-2-2c)$. Observe that
\[
\log\frac{\omega^n_{\varphi}}{\omega^n}\omega^n_{\varphi}\geq-e^{-1}\omega^n,
\]
We have \begin{eqnarray*} c&&\geq
\left(\int_{\epsilon\omega^n_{\varphi}>\omega^n}+\int_{\epsilon\omega^n_{\varphi}\leq\omega^n}\right)
\left(\log\frac{\omega^n_{\varphi}}{\omega^n}\omega^n_{\varphi}\right)\\
&&\geq\int_{\epsilon\omega^n_{\varphi}>\omega^n}\left(\log\frac{1}{\epsilon}\right)\omega^n_{\varphi}
+\int_{\epsilon\omega^n_{\varphi}\leq\omega^n}\left(-e^{-1}\omega^n\right)\\
&&>2(1+c)\int_{\epsilon\omega^n_{\varphi}>\omega^n}\omega^n_{\varphi}-1.\end{eqnarray*}
It follows that
\[\int_{\epsilon\omega^n_{\varphi}>\omega^n}\o^n_\varphi<\frac{1}{2}.\]
And so
\begin{equation}\label{5-6}\int_{\epsilon\omega^n_{\varphi}\leq\omega^n}\o^n_\varphi>\frac{1}{2}.\end{equation}
Also we have
\begin{eqnarray}\label{5-7}\int_{\omega^n\leq4\omega^n_{\varphi}}\omega^n_{\varphi}
&=&1-\int_{\o^n
>4\o^n_{\varphi}}\o^n_\varphi\nonumber\\&>&1-\int_{\o^n>4\o^n_{\varphi}}\frac{\o^n}{4}\geq
1-\frac{1}{4}\int_M\o^n= \frac{3}{4}.
\end{eqnarray}
And we know \begin{equation}\label{5-8}
\int_{\omega^n\leq4\omega^n_{\varphi}}\omega^n\geq\int_{\frac{\epsilon}{4}\omega^n\leq
\epsilon\omega^n_{\varphi}\leq\omega^n}\o^n\geq\epsilon\int_{\frac{\epsilon}{4}\omega^n\leq
\epsilon\omega^n_{\varphi}\leq\omega^n}\omega^n_{\varphi}.\end{equation}

Set $A=\{x\in M: \epsilon\omega^n_{\varphi}(x)\leq\omega^n(x)\}$
and $B=\{x\in M: \omega^n(x)\leq4\omega^n_{\varphi}(x)\}$. By
(\ref{5-6}) (\ref{5-7}),
\begin{equation}\label{5-9}
\int_{\frac{\epsilon}{4}\omega^n\leq
\epsilon\omega^n_{\varphi}\leq\omega^n}\omega^n_{\varphi}=\int_{A\cap
B}\o_\varphi^n=\int_A\o_\varphi^n+\int_B\o_\varphi^n-\int_{A\cup
B}\o^n_\varphi>\frac{1}{2}+\frac{3}{4}-1=\frac{1}{4}.
\end{equation}
Combine (\ref{5-8}) and (\ref{5-9}), we have
\begin{equation}\label{5-10}\int_{\omega^n\leq4\omega^n_{\varphi}}\omega^n>
\epsilon/4.
\end{equation}
By Green's formula, we have \begin{equation}\label{5-11}
F(p)=-\int_MG(p,q)\triangle
F(q)\omega^n(q)+\int_MF\omega^n,\end{equation} where $G(p,q)\geq
0$ is a Green function of $\omega$. By (\ref{5-10}) (\ref{5-11})
and the uniform bound of $\triangle F$,
\begin{eqnarray*}\inf_MF&&\geq\int_M F\omega^n-C\\
&&\geq\inf_MF\int_{\omega^n\geq4\omega^n_{\varphi}}\omega^n+\int_{\omega^n<4\omega^n_{\varphi}}F\omega^n-C\\
&&\geq\inf_MF\int_{\omega^n\geq4\omega^n_{\varphi}}\omega^n-\log4\int_{\omega^n<4\omega^n_{\varphi}}\omega^n-C\\
&&\geq\left(1-\frac{\epsilon}{4}\right)\inf_MF-C,\end{eqnarray*}
where we can assume $\inf_M F<0.$ Therefore, we have
\[\inf_MF\geq-4C\exp(2+2c).\] By the way we choose the constant
$c$ in the beginning of the proof, the proposition is proved.
\end{proof}
So we get a uniform bound for $F$, $\triangle F$ and $\triangle
\varphi$. It implies that all the metrics $g_\varphi$ are
equivalent. For the Monge-Ampere equation
\[\frac{\det\left(g_{i\bar{j}}+\varphi,_{i\bar{j}}\right)}{\det\left(g_{i\bar{j}}\right)}=\exp(F),\]
The $C^2$ bound of $F$ implies that $\varphi \in W^{4,p}(M)$ is
uniformly bounded for any $p>1.$  By the Sobolev's embedding
theorem, for any $\alpha \in(0,1)$, $\varphi \in C^{3,\alpha}$ is
uniformly bounded.
\end{proof}

\begin{cor}For the Calabi flow initiating from any smooth K\"ahler metric,
the flow exists as long as the Ricci curvature stays uniformly
bounded.
\end{cor}
\begin{proof}
By the short time existence of the Calabi flow, we can assume the
maximal existence interval of the flow is $[0, T)$. If $T<\infty$,
because $Ric_{\varphi(t)}$ is uniformly bounded, so is the scalar
curvature. Then by \begin{eqnarray*} \frac{\p \varphi}{\p
t}=R_\varphi-\underline{R},\end{eqnarray*} we know $\varphi$ is
uniformly bounded. By Theorem 5.1, we know $\varphi(t)$ is
uniformly bounded in $C^{3, \alpha}$ for any $\alpha \in (0,1)$
and all metrics along the flow are equivalent. By Theorem 3.2, the
flow can be extended. Contradiction.
\end{proof}
\begin{cor}
The maximal existence interval of the Calabi flow is  $[0, T)$,
starting with any initial smooth potential $\varphi(0)$. If the
Ricci curvature and the potential are both uniformly bounded in
$[0, T)$, then $T=\infty$ and the flow converges to some extremal
metric in the same K\"ahler class. Moreover, when the limit metric
is a CscK metric, the convergence is exponentially fast.
\end{cor}
\begin{proof}Obviously $T=\infty$. To get convergence, first we show the flow converges by sequence.
Because Ricci curvature $R_{i\bar{j}}(t)$ and the potential
$\varphi(t)$ are both uniformly bounded, all metrics $g(t)$ are
equivalent and $\varphi(t)$ is uniformly bounded in $C^{3,
\alpha}(M)$ for any $t \in [0, \infty)$ and any $0<\alpha<1$. Then
by the smoothing property of the Calabi flow (Theorem 3.3), we can
get actually $\varphi(t)$ is uniformly bounded in $C^{k,
\alpha}(M)$ for any $k \in \mathbb{N}$. It follows that the flow
converges by sequence in $C^{\infty}$. To get the limit metric is
extremal metric, recall that the Calabi energy is decreasing along
the flow, \begin{eqnarray*}
\frac{d}{dt}\int_M(R(t)-\underline{R})^2dg(t)=-2\int_MR(t),_{\alpha\beta}R(t)^{,\alpha\beta}dg(t).\end{eqnarray*}
Denote \begin{eqnarray*} E_0=\inf_{[0,
\infty)}\int_M(R(t)-\underline{R})^2dg(t).\end{eqnarray*} For any
sequence $t_n\rightarrow \infty$, we have a subsequence $t_{n_k}$,
such that $\varphi(t_{n_k})$ converges to $\varphi_{\infty}$ in
$C^{\infty}$ topology. Then we have \begin{eqnarray*} E_0=
\int_M(R_\infty-\underline{R})^2dg_\infty.\end{eqnarray*}
Moreover,
\begin{eqnarray*} E_{t_n}-E_0&=&\lim_{k\rightarrow
\infty}\left(\int_M(R(t_n)-\underline{R})^2dg(t_n)-\int_M(R(t_{n_k})-\underline{R})^2dg(t_{n_k})\right)\\
&=&\lim_{k\rightarrow \infty}\int_{t_n}^{t_{n_k}}\frac{d}{dt}
\int_M(R(t)-\underline{R})^2dg(t)dt\\
&=&-2\lim_{k\rightarrow
\infty}\int_{t_n}^{t_{n_k}}\int_MR(t),_{\alpha\beta}R(t)^{,\alpha\beta}dg(t)dt\\
&=&-2\int_{t_n}^{\infty}\int_MR(t),_{\alpha\beta}R(t)^{,\alpha\beta}dg(t)dt.
\end{eqnarray*} When $t_n\rightarrow \infty,$ $E_{t_n}\rightarrow E_0$, it
implies $
\int_MR(t),_{\alpha\beta}R(t)^{,\alpha\beta}dg(t)\rightarrow 0$
when $t\rightarrow 0.$ In particular,
\begin{eqnarray*}\int_MR_\infty,_{\alpha\beta}{R_\infty}^{,\alpha\beta}dg_\infty=0.\end{eqnarray*}
So $R_\infty,_{\alpha\beta}=0$, $g_\infty$ is extremal metric.
When $R_\infty=\underline{R},$ by similar argument as in Theorem
4.1, the convergence is exponentially fast.
\end{proof}
\section{Removing singularity}
In this section we consider a K\"ahler metric which is defined
locally in the punctured disc $D\backslash\{0\}$, where $D=\{z \in
\mathbb{C}^n: |z|\leq 1\}$. If the metric admits constant scalar
curvature in some weak sense, we show the metric can be extended
to the whole disc smoothly.
\begin{defi}A K\"ahler metric
$g_\varphi=(\delta_{ij}+\p_i\p_{\bar{j}}\varphi)dz^i\otimes
dz^{\bar{j}}$ is defined in the punctured disc $D\backslash
\{0\}$, if $\varphi \in C^{2}(D\backslash \{0\})$ is a real valued
function such that $\{\delta_{ij}+\p_i\p_{\bar{j}}\varphi\}$ is a
positive Hermitian matrix function. The scalar curvature
$R_\varphi$ of $g_\varphi$ is constant in weak sense in the
punctured disc if $\varphi \in C^{2}(D\backslash\{0\})\cap W^{3,
2}_{loc}(D\backslash \{0\})$ satisfies \begin{eqnarray*}
-g^{i\bar{j}}_{\varphi}\p_i\p_{\bar{j}}\log{\det
(\delta_{kl}+\p_k\p_{\bar{l}}\varphi)}=\underline{R},\quad
\mbox{in}\quad D\backslash \{0\}\end{eqnarray*} in weak sense.
\end{defi}
\begin{theo}
If $g_\varphi$ is a K\"ahler metric defined in the punctured disc
with constant scalar curvature in weak sense in the punctured
disc, and $\p_i\p_{\bar{j}} \varphi \in L^{\infty}(D)$,
$L=g^{i\bar{j}}_{\varphi}\p_i\p_{\bar{j}}$ is uniformly elliptic
on $D$, moreover, to avoid regularity issues on the boundary, we
assume $\varphi$ is smooth in the neighborhood of the boundary $\p
D$, then $g_\varphi$ can be extended to a smooth metric with
constant scalar curvature in the whole disc $D$.
\end{theo}
We need to prove a theorem first:
\begin{theo} Let $ \left(a^{i\b j}\right)$ be a uniformly elliptic matrix in  $D$.
The following Dirichlet problem can be solved in $W^{1, 2}(D)$:
 \begin{eqnarray*}
-\p_i\left(a^{i\bar{j}} \det (a_{k\b l})  \p_{\b j} u\right)&=&
\underline{R}\det{(a_{i\b j})}~~\mbox{in}~~ D,\\
u|_{\p D}&=&h.\end{eqnarray*} Here $\left(a^{i\b j}\right) \cdot
\left(a_{k\b l}\right) = I_{n \times n}.$
\end{theo}
\begin{proof}To solve the Dirichlet equation, define
the following functional: \begin{eqnarray*}
I(u)=\int_D\left(a^{i\bar{j}}\p_iu\p_{\bar{j}}u+\underline{R}
u\right)\det{(a_{i\bar{j}})}dx,\end{eqnarray*} where $u \in
A=\left\{u|u \in W^{1,2}(D), u|_{\partial D}=h\right\}$, and $dx$
is the Euclidean measure. Here $u|_{\partial D}=h$ is in trace
sense \cite{Ev}.
\begin{claim}$I(u)$ is bounded from below.\end{claim}
\begin{proof} $\forall u \in A$ \begin{eqnarray*}
I(u)&=&\int_D\left(a^{i\bar{j}}\p_iu\p_{\bar{j}}u+\underline{R} u\right)\det{(a_{i\bar{j}})}dx\\
&\geq&\frac{1}{C}\int_D|Du|^2dx-C\int_D |u|dx\\
&\geq&\frac{1}{C}\int_D|Du|^2dx-\epsilon
C\int_Du^2dx-\frac{C}{\epsilon},\end{eqnarray*} where $D$ is the
derivative under the Euclidean metric. Let $\o\in A$, we know
\begin{eqnarray*}
\int_D|Du|^2dx&\geq&\int_D|Du-Dw|^2dx-\int_D|Dw|^2dx,\\
\int_Du^2dx&\leq&2\int_D(u-w)^2dx+2\int_Dw^2dx.\end{eqnarray*} So
we have
\begin{eqnarray*}
I(u)&\geq&\frac{1}{C}\int_D\left(|D(u-w)|^2-|Dw|^2\right)dx-2\epsilon
C\int_D(u-w)^2dx-2\epsilon C\int_Dw^2dx-\frac{C}{\epsilon}\\
&\geq&\frac{1}{C}\int_D|D(u-w)|^2dx-2\epsilon
C\int_D(u-w)^2dx-C(\epsilon, w)\\
&\geq&\frac{1}{2C}\int_D|D(u-w)|^2dx-C(\epsilon, w)
\\ &\geq&-C(\epsilon,w).\end{eqnarray*} The above inequality follows from
Poincar\'e inequality because $u-w \in W^{1,2}_0(D)$ if
$\epsilon>0$ small enough. Denote $m=\inf_{u \in A}I(u)$. Taking a
minimizing sequence $\{u_k\}$ of $I(u)$, we have \begin{eqnarray*}
\frac{1}{2C}\int_D|D(u_k-w)|^2dx\leq I(u_k)+C(\epsilon,
w),\end{eqnarray*} for some fixed $w \in A$. It is easy to show
that $u_k$ is uniformly bounded in $W^{1,2}(D)$. So we can get a
subsequence of $\{u_k\}$ converge to $u_0$ weakly in $W^{1,2}(D)$,
and strongly in $L^2(D)$.
\end{proof}

\begin{claim}
$u_0$ is the minimizer of $I(u)$ and $u_k$ converge to $u_0$
strongly in $W^{1,2}(D)$.\end{claim}
\begin{proof}
We know
\begin{eqnarray*}\int_Da^{i\bar{j}}\p_i(u_k-u_0)\p_{\bar{j}}(u_k-u_0)\det{(a_{i\bar{j}})}dx\geq0.\end{eqnarray*}
It implies
\begin{eqnarray*}\liminf_k\int_Da^{i\bar{j}}\p_i(u_k)\p_{\bar{j}}(u_k)\det{(a_{i\bar{j}})}dx\geq
\int_Da^{i\bar{j}}\p_i(u_0)\p_{\bar{j}}(u_0)\det{(a_{i\bar{j}})}dx,
\end{eqnarray*} because $u_k$ converges to $u_0$ weakly in $W^{1,2}(D).$ Then
we
have \begin{eqnarray*} m&=&\lim_kI(u_k)\\
&=&\lim_k\int_D\left(a^{i\bar{j}}\p_iu_k\p_{\bar{j}}u_k+\underline{R}
u_k\right)\det{(a_{i\bar{j}})}dx\\
&\geq&\int_D\left(a^{i\bar{j}}\p_iu_0\p_{\bar{j}}u_0+\underline{R}
u_0\right)\det{(a_{i\bar{j}})}dx\\
&=&I(u_0).\end{eqnarray*} $m$ is the minimum of $I(u)$, so $u_0$
is the minimizer. It easily follows that $u_k$ converge to $u_0$
strongly.
\end{proof}
The minimizer $u_0$ of $I(u)$ satisfies the following equation in
weak sense: \begin{eqnarray*}
-\p_i\left(a^{i\bar{j}}\det{(a_{k\bar{l}})}\p_{\bar{j}}u\right)&=&
\underline{R}\det{(a_{i\bar{j}})}~~\mbox{in}~~ D.\end{eqnarray*}
Theorem 6.3 follows.\end{proof}

Now we are in the position to prove Theorem 6.2.

\begin{proof}
Taking $a_{i\bar{j}}=\delta_{ij}+\p_i\p_{\bar{j}}\varphi,$
$h=\log{\det{(\delta_{ij}+\p_i\p_{\bar{j}}\varphi)}}|_{\p D}$,
Theorem 6.3 says there exists a weak solution $u_0$ of the
following equation,
\begin{eqnarray*}
-\p_i\left(g^{i\bar{j}}_{\varphi}\det{(\delta_{kl}+\p_k\p_{\bar{l}}\varphi)}\p_{\bar{j}}u\right)&=&
\underline{R}\det{(\delta_{ij}+\p_i\p_{\bar{j}}\varphi)}~~\mbox{in}~~
D,\\u|_{\p D}&=&h.\end{eqnarray*} Since $\varphi \in
C^{2}(D\backslash\{0\})\cap W^{3, 2}_{loc}(D\backslash \{0\})$ ,
by straightforward calculation, $u_0$ satisfies the following
equation in weak sense:\begin{eqnarray*}
-g^{i\bar{j}}_{\varphi}\p_i\p_{\bar{j}}u_0=\underline{R}~~\mbox{in}~~
D.\end{eqnarray*} Denote
$v=\log{\det{(\delta_{ij}+\p_i\p_{\bar{j}}\varphi)}}-u_0$, where
$v$ is the weak solution of
\begin{eqnarray*}-g^{i\bar{j}}_{\varphi}\p_i\p_{\bar{j}}v&=&0,~~\mbox{in}~D\backslash\{0\},\\
v|_{\partial D}&=&0.\end{eqnarray*} If we can show $v$ is
identically zero, then $\varphi$ will satisfy the following
equation
\begin{equation} \label{6-1}-g^{i\bar{j}}_{\varphi}\p_i\p_{\bar{j}}\log{\det
(\delta_{kl}+\p_k\p_{\bar{l}}\varphi)}=\underline{R}\end{equation}
in the whole disc. It implies that we can extend the metric
$g_\varphi$ to the whole
disc with constant scalar curvature in weak sense.\\

To show $v$ is actually 0, we need to use maximum principle for
elliptic equation with an isolated singularity \cite{Gs}. In
general, we need some further assumption for the elliptic operator
around the singularity (cf. \cite{Gs}). In this complex setting,
it is nice that we don't need any assumption
other than the uniform ellipticity.\\

Let $f= |z|^p$, then
\[
 f_i = p |z|_i |z|^{p-1}  = {p \over 2} \b {z_i} |z|^{p-2}.
\]
and
\[
f_{i \b j} = {p \over 2} |z|^{p-2} \delta_{i j} + { p(p-2)\over 4}
\b {z_i} z_j |z|^{p-4}.
\]
Set
\[
\psi = v f, \qquad {\rm or}\;\;v = \psi |z|^{-p} = \psi
|z|^q,\qquad q=-p < 0.
\]
Then,
\[v_i = \psi_i |z|^q + \psi  {q \over 2} \b {z_i} |z|^{q-2}
\]
and
\[
v_{i \b j} = \psi_{i \b j} |z|^q + \psi_{\b j}  {q \over 2} \b
{z_i} |z|^{q-2} + \psi _i  {q \over 2}  {z_j} |z|^{q-2} + \psi
\left(  {q \over 2} |z|^{q-2} \delta_{i j} + { q(q-2)\over 4} \b
{z_i} z_j |z|^{q-4}\right).
\]
Claim that $\psi \leq 0.\;$ Otherwise, there is at least one
interior point $|z|\neq 0$ such that
\[
\left(\psi_{i \b j}\right) \leq 0, \qquad {\rm and}\;\; \psi > 0.
\]
Note that $g^{i \b j} v_{i \b j} = 0.\;$ Thus, we have (at maximum
point)
\[
0 \leq g^{i\b j} \psi \left(  {q \over 2} |z|^{q-2} \delta_{i j} +
{ q(q-2)\over 4} \b {z_i} z_j |z|^{q-4}\right).
\]
Since $q< 0$ and $|z| \neq 0$, this means that ($q < 0$)
\[
 |z|^{2} \sum_{i=1}^n g^{i \b i} + { (q-2)\over 2} g^{ i\b j}  \b {z_i} z_j  \leq 0.
\]
This is a contradiction  when $-q > 0$ is very small.   Similarly,
we can show that $\psi \geq 0.$ Consequently, $\psi = 0.\;$
\begin{claim}
$g_{\varphi}$ is a smooth metric with constant scalar curvature in
the entire disc.
\end{claim}
\begin{proof} By (\ref{6-1}), $\varphi$ is the weak solution of
\[
-g^{i\bar{j}}_{\varphi}\p_i\p_{\bar{j}}\log{\det
(\delta_{kl}+\p_k\p_{\bar{l}}\varphi)}=\underline{R}.
\]
Since $L=g^{i\bar{j}}_{\varphi}\p_i\p_{\bar{j}}$ is uniformly
elliptic, the De Giorgi-Nash-Moser's estimate \cite{Gt} gives that
\begin{equation}\label{6-2}
u=\log \det{(\delta_{i\bar j}+\varphi_{i\bar j})}
\end{equation}
is in $C^{\alpha}(D)$. Also by the Evans-Krylov estimate
\cite{Evans01} for Monge-Ampere equation
($L=g^{i\bar{j}}_{\varphi}\p_i\p_{\bar{j}}$ is uniformly
elliptic), we can deduce from (\ref{6-2}) that $\varphi\in C^{2,
\alpha}(D)$. Go back to (\ref{6-1}),
\[
-g^{i\bar{j}}_{\varphi}\p_i\p_{\bar{j}}\log{\det
(\delta_{kl}+\p_k\p_{\bar{l}}\varphi)}=\underline{R}.
\]
Now $g_\varphi\in C^{\alpha}(D)$, and so the elliptic theory gives
that $u=\log{\det{(\delta_{i\bar j}+\varphi_{i\bar j})}}\in C^{2,
\alpha}(D)$. The standard Monge-Ampere theory implies that
$\varphi\in C^{4, \alpha}$. Using  the standard boot-strapping
argument, we can get that $\varphi\in C^{\infty}(D)$. To avoid the
issue of boundary regularity, we assume in addition the $\varphi$
is smooth in the neighborhood of the boundary $\p D$. So
$\varphi\in C^{\infty}(\bar D).$
\end{proof}
\end{proof}

\end{document}